\documentclass[10pt]{article}

\usepackage{amsmath} % eg. for \begin{equation}
\usepackage[mathscr]{eucal}
\usepackage{amssymb} % eg. for \leqslant
\usepackage{theorem}
\usepackage{tikz-cd}

\theoremstyle{change}  % puts numbers IN FRONT of "Theorem"
\newtheorem{theorem}{Theorem}[section] % defines environment "Theorem".
\newtheorem{lemma}[theorem]{Lemma}  % defines environment "Lemma", that
\newtheorem{proposition}[theorem]{Proposition}
\newtheorem{corollary}[theorem]{Corollary}
\theorembodyfont{\rmfamily}  % has the effect that the content of Remark,
\newtheorem{Remark}[theorem]{Remark}

\newtheorem{definition}[theorem]{Definition}

\newtheorem{nothing}[theorem]{} % empty Theoremumgebung.

\newenvironment{proof}{\noindent{\bf Proof}\ }{\qed\bigskip}

\renewcommand{\le}{\leqslant} % needs amssymb-Paket

\newcommand{\btilde}{\tilde{b}}

\newcommand{\catfont}{\mathsf}

\newcommand{\ctilde}{\tilde{c}}

\newcommand{\End}{\mathrm{End}}

\newcommand{\etilde}{\tilde{e}}

\newcommand{\FF}{\mathbb{F}}

\newcommand{\ftilde}{\tilde{f}}
\newcommand{\Gal}{\mathrm{Gal}}

\newcommand{\Hom}{\mathrm{Hom}}

\newcommand{\Ind}{\mathrm{Ind}}

\newcommand{\lexp}[2]{\setbox0=\hbox{$#2$} \setbox1=\vbox to
                 \ht0{}\,\box1^{#1}\!#2}

\newcommand{\lmod}[1]{\llap{\phantom{|}}_{#1}\catfont{mod}}

\newcommand{\myiso}{\buildrel\sim\over\to}

\newcommand{\qed}{\nobreak\hfill
                  \vbox{\hrule\hbox{\vrule\hbox to 5pt
                  {\vbox to 8pt{\vfil}\hfil}\vrule}\hrule}}
\newcommand{\Res}{\mathrm{Res}}

\newcommand{\stab}{\mathrm{stab}}
\newcommand{\stilde}{\tilde{s}}
\newcommand{\Stilde}{\tilde{S}}

\newcommand{\Ttilde}{\tilde{T}}

\newcommand{\tr}{\mathrm{tr}}

\newcommand{\Vtilde}{\widetilde{V}}

\newcommand{\Wtilde}{\widetilde{W}}

\newcommand{\ZZ}{\mathbb{Z}}

\newcommand{\Tr}{\mathrm{Tr}}
\newcommand{\Br}{\mathrm{Br}}

\title{Galois descent of equivalences between blocks of $p$-nilpotent groups
\footnote{{\bf MR Subject Classification:}  
20C20, 19A22. {\bf Keywords:}  $p$-permutation modules, trivial source modules, splendid Rickard equivalence, $p$-permutation equivalence, $p$-nilpotent groups, Galois descent.}}
\author{\small Robert Boltje\\
  \small Department of Mathematics\\
   \small University of California\\
   \small Santa Cruz, CA 95064\\
   \small U.S.A.\\
   \small boltje@ucsc.edu
   \and 
   \small Deniz Y\i lmaz\\
   \small Institut f\"ur Mathematik\\
   \small Friedrich-Schiller-Universit\"at Jena\\
   \small 07737 Jena\\
   \small Germany\\
   \small deyilmaz@ucsc.edu}
\date{February 25, 2021 (revised May 1, 2021)}

\begin{document}

\sloppy

\maketitle

%%%%%%%%%%%%%%%%%%%%%%%%%%%% ABSTRACT %%%%%%%%%%%%%%%%%%%%%%%%%%%%%%%%

\begin{abstract}
We give sufficient conditions on $p$-blocks of $p$-nilpotent groups over $\FF_p$ to be splendidly Rickard equivalent and $p$-permutation equivalent to their Brauer correspondents. The paper also contains Galois descent results on $p$-permutation modules and $p$-permutation equivalences that hold for arbitrary groups.
\end{abstract}

%%%%%%%%%%%%%%%%%%% SECTION 1 %%%%%%%%%%%%%%%%%%%%%%%%%%%%%%%%%%%%%%

\section{Introduction}

In \cite{KL2018}, Kessar and Linckelmann proved that Brou\'e's Abelian Defect Group Conjecture (originally stated over splitting fields) holds for blocks with cyclic defect groups over {\em arbitrary} fields of characteristic $p>0$, in particular over the prime field $\FF_p$. More precisely, if $G$ is a finite group and $b$ is a block idempotent of $\FF_p G$ with cyclic defect group $D$, then there exists a {\em splendid Rickard equivalence} between $\FF_p Gb$ and its Brauer correspondent block algebra $\FF_pH\Br_D(b)$, where $H=N_G(D)$ and $\Br_D\colon (\FF_pG)^D\to\FF_pC_G(D)$ is the Brauer homomorphism, an $\FF_p$-algebra homomorphism which is given by truncation. 

\medskip
In this paper we investigate if a similar phenomenon holds for blocks of $p$-nilpotent groups. In this case, a positive answer over a splitting field $F$ of characteristic $p>0$ of $G$ was given by Rickard, see \cite{Rickard1996}, even without the assumption of abelian defect groups. There, he introduced and used the notion of an {\em endosplit $p$-permutation resolution} in order to construct such splendid Rickard equivalences. We have two main results. The first gives sufficient conditions under which there exists such a splendid Rickard equivalence between Brauer corresponding blocks of a $p$-nilpotent group over $\FF_p$. The second gives sufficient conditions under which the weaker form of equivalence, namely a {\em $p$-permutation equivalence}, exists over $\FF_p$.

\medskip
So let $G$ be a {\em $p$-nilpotent group}, i.e., a finite group whose largest normal $p'$-subgroup $N$ is a complement to a (and then each) Sylow $p$-subgroup of $G$. Moreover, let $F$ be a finite splitting field of $G$ and its subgroups of characteristic $p>0$. Let $\btilde$ be a block idempotent of $\FF_pG$ and let $b$ be a block idempotent of $FG$ which occurs in a primitive decomposition of $\btilde$ in $Z(FG)$. Then $b$ is contained in $Z(FN)$. Let $e$ be a block idempotent of $FN$ that occurs in the primitive decomposition of $b$ in $Z(FN)$. Adjoining the coefficients of $e$ and of $b$ to $\FF_p$, one obtains subfields $\FF_p[b]\subseteq \FF_p[e]\subseteq F$, since $b$ is a sum of $G$-conjugates of $e$.

\bigskip\noindent
{\bf Theorem A}\quad {\it Let $G$ be a $p$-nilpotent group and let $\btilde$ be a block idempotent of $\FF_pG$. Suppose that $p$ is odd or $\btilde$ has abelian defect groups, and suppose that, with the above notation, $\FF_p[b]=\FF_p[e]$. Then there exists a splendid Rickard equivalence between the block algebra $\FF_pG\btilde$ and its Brauer correspondent block algebra.}

\bigskip
Theorem~A follows from the more precise statement in Proposition~\ref{prop Rickard equivalence over Fp} and Remark~\ref{rem on hypo A}. The proof uses Rickard's original approach in \cite{Rickard1996} involving endosplit $p$-permutation resolutions, a descent result in \cite{KL2018}, and the classification of endopermutation modules over $p$-groups, see \cite{Thevenaz2007} for a survey article on the latter.

\bigskip
There are weaker forms of equivalences between blocks than splendid Rickard equivalences, as for instance $p$-permutation equivalences which were introduced in \cite{BX2008} and extended in \cite{Linckelmann2009} and \cite{BoltjePerepelitsky2020}. See Section~\ref{sec descent of equivalences} for a definition.

%
%\medskip
%A weaker form of Brou\'e's Abelian Defect Group Conjecture states that if $b$ is a block idempotent of a group algebra $FG$ with abelian defect groups then there exists a {\em $p$-permutation equivalence} between the block algebra $FGb$ and its Brauer correspondent block algebra. The notion of a $p$-permutation equivalence was introduced in \cite{BX2008} and extended in \cite{Linckelmann2009} and \cite{BoltjePerepelitsky2020}.

\bigskip\noindent
{\bf Theorem B}\quad {\it Let $G$ be a $p$-nilpotent group with abelian Sylow $p$-subgroup and let $\btilde$ be a block idempotent of $\FF_p G$. Then there exists a $p$-permutation equivalence between $\FF_pG\btilde$ and its Brauer correspondent block algebra.}

\bigskip
Theorem B follows from the more precise statement in Corollary~\ref{cor p-perm equiv}. 
The proof uses again Rickard's construction and Galois descent arguments for the representation ring of trivial source modules developed in this paper, see Theorem~\ref{thm Galois descent} and Lemma~\ref{lem descending}. The reason that we only obtain a $p$-permutation equivalence and not a splendid Rickard equivalence in Theorem B, is that we don't have a descent result analogous to Lemma~\ref{lem descending} for splendid Rickard equivalences and that the descent result from \cite{KL2018} cannot be applied without the assumption that $\FF_p[b]=\FF_p[e]$, see also Remark~\ref{rem no descent theory for Z}. Because Theorem~\ref{thm Galois descent} and Lemma~\ref{lem descending} are of independent interest we include them in the introduction as Theorems~C and D. For these two results, $G$ and $H$ can be arbitrary finite groups and we assume that $F$ is a splitting field for $G$ and $H$ and their subgroups.

\bigskip\noindent
{\bf Theorem C}\quad{\it
Let $F'$ be a subfield of $F$ and $\Delta:=\Gal(F/F')$. Then, scalar extension from $F'$ to $F$ induces an isomorphism $T_{F'}(G)\to T_F(G)^\Delta$ from the trivial source ring of $F'G$ to the $\Delta$-fixed points of $T_F(G)$.}

\bigskip\noindent
{\bf Theorem D}\quad{\it
Let $b$ and $c$ be block idempotents of $FG$ and $FH$, respectively. Let $\btilde$ and $\ctilde$ denote the block idempotents of $\FF_pG$ and $\FF_pH$ associated to $b$ and $c$, respectively, as in Proposition~\ref{prop Gamma and blocks}(a). Moreover, let $\omega\in T^\Delta(FGb,FHc)$ be a $p$-permutation equivalence between $FGb$ and $FHc$. Suppose that $\stab_{\Gamma}(\omega)=\stab_{\Gamma}(b)=\stab_{\Gamma}(c)$. Then there exists a $p$-permutation equivalence between $\FF_pG\btilde$ and $\FF_pH\ctilde$.
}

\bigskip
The paper is arranged as follows. In Section~2 we prove Theorem~\ref{thm Galois descent}. The definition and basic properties of endosplit $p$-permutation resolutions are given in Section~3. In Section~4 we collect basic results on the Galois group action on blocks and prove Lemma~\ref{lem descending}. Finally, in Section~5 we prove Theorems A and B.

\medskip
Our notation is standard. For any rings $R$ and $S$ we denote by $\lmod{R}$ (resp.~$\lmod{R}_S$) the categories of finitely generated left $R$-modules (resp.~$(R,S)$-bimodules). For objects $M$ and $N$ in a module category or chain complex category we write $M\mid N$ to indicate that $M$ is isomorphic to a direct summand of $N$. If $H$ and $K$ are subgroups of a finite group $G$, then $g\in G/H$ (resp.~$g\in H\setminus G/K$) indicates that $g$ runs through a set of representatives of the given cosets (resp.~double cosets).

\bigskip\noindent
{\bf Acknowledgement} The authors are most grateful to the referee for her/his thorough reading and detailed comments which among other improvements resulted in a correction of the statement of Proposition~\ref{prop stableunderT}.

%%%%%%%%%%%%%%%%%%%% SECTION 2 %%%%%%%%%%%%%%%%%%%%%%%%%%%%%%%%%%%%%

\section{Galois descent of $p$-permutation modules}\label{sec descent of modules}

Throughout this paper, $G$ and $H$ denote finite groups and $F$ a finite field of characteristic $p$ which is a splitting field for all subgroups of $H$ and $G$. Moreover, $\Gamma:=\Gal(F/\FF_p)$ denotes the Galois group of $F$ over $\FF_p$. For any subfield $F'\subseteq F$ one has functors 
\begin{equation*}
   -_{F'}\colon \lmod{FG}\to \lmod{F'G}\quad\text{and}\quad F\otimes_{F'}-\colon \lmod{F'G}\to \lmod{FG}
\end{equation*}
defined by restriction and extension of scalars.

\begin{proposition}\label{prop vertices} Let $Q$ be a $p$-subgroup of $G$ and let $F'$ be a subfield of $F$. If $M\in \lmod{F'G}$ is relatively $Q$-projective then $F\otimes_{F'}M\in\lmod{FG}$ is relatively $Q$-projective.
If $N\in\lmod{FG}$ is relatively $Q$-projective then $N_{F'}\in\lmod{F'G}$ is relatively $Q$-projective. 
\end{proposition}

\begin{proof}
This follows immediately from the fact that restriction and extension of scalars commute with $\Ind_Q^G$.
\end{proof}

For each $\sigma\in \Gamma$ one has a functor
\begin{equation}\label{eqn Galois functors}
   \lexp{\sigma}{-} \colon \lmod{FG}\to\lmod{FG}
\end{equation}
which assigns to $M\in\lmod{FG}$ the $FG$-module $\lexp{\sigma}{M}$ whose underlying abelian group is equal to $M$ and whose $FG$-module structure is given by restriction along the ring isomorphism $\sigma^{-1}\colon FG\to FG$, $\sum_{g\in G}\alpha_g g\mapsto \sum_{g\in G} \sigma^{-1}(\alpha_g) g$. For any $FG$-module homomorphism $f$ one has $\lexp{\sigma}{f}=f$. Similarly one defines the functor $\lexp{\sigma}{-}\colon \lmod{FG}_{FH}\to \lmod{FG}_{FH}$. For any $M\in\lmod{FG}$ we set $\stab_\Gamma(M):=\{\sigma\in \Gamma\mid \lexp{\sigma}{M}\cong M\}$.

We recall from \cite[Definition~5.4.10]{Linckelmann2018a} the definition of the Brauer construction functor 
\begin{equation*}
  -(P)\colon \lmod{FG}\to\lmod{F[N_G(P)/P]}\,,
\end{equation*}
for any $p$-subgroup $P$ of $G$, and we denote by
\begin{equation*}
  -^\circ:=\Hom_F(-,F)\colon \lmod{FG}\to\lmod{FG}
\end{equation*} 
the functor of taking $F$-duals. The above functors extend to functors between appropriate categories of (co-)chain complexes and they have the following properties.

\begin{lemma}\label{lem functor properties}
Let $G$, $H$ and $K$ be finite groups. Further, let $M$ and $N$ be $FG$-modules, $U$ an $(FG,FH)$-bimodule, $V$ an $(FH,FK)$-bimodule, $L\le G$ a subgroup, $P\le G$ a $p$-subgroup and $\sigma,\tau \in \Gamma$. Moreover, let $F'\subseteq F$ be a subfield and set $\Delta:=\Gal(F/F')$. Then one has

\smallskip
{\rm (a)} $\lexp{\tau\circ\sigma}M = \lexp{\tau}{\left(\lexp{\sigma}M\right)}$,  $\lexp{\sigma}{\left(M\oplus N\right)} = \lexp{\sigma}M\oplus \lexp{\sigma}N$, and $\lexp{\sigma}{\left(M\otimes_{F}N\right)} = \lexp{\sigma}M\otimes_{F} \lexp{\sigma}N$.

\smallskip
{\rm (b)} $\Res_L^G(\lexp{\sigma}M)=\lexp{\sigma}{\left(\Res_L^G M\right)}$ and $\Ind_L^G(\lexp{\sigma}M)=\lexp{\sigma}{\left(\Ind_L^G M\right)}$.

\smallskip
{\rm (c)} $\left(\lexp{\sigma}M\right)^\circ = \lexp{\sigma}{\left(M^\circ\right)}$, $\lexp{\sigma}{\left(U\otimes_{FH}V\right)} = \lexp{\sigma}U\otimes_{FH} \lexp{\sigma}V$, and $\left(\lexp{\sigma}M\right)(P)=\lexp{\sigma}{\left((M)(P)\right)}$.

\smallskip
{\rm (d)} $F\otimes_{F'} M_{F'}\cong \bigoplus_{\sigma\in \Delta} \lexp{\sigma}M$.
\end{lemma}

\begin{proof}
The proofs of (a)--(c) are straightforward. For a proof of Part~(d) see \cite[Proposition 6.3]{KL2018}.
\end{proof}

\begin{corollary}\label{cor vertex}
Let $F'$ be a subfield of $F$. An indecomposable $F'G$-module $M$ and the indecomposable direct summands of $F\otimes_{F'G} M$ have the same vertices. Similarly, an indecomposable $FG$-module $N$ and the indecomposable direct summands of $N_{F'}$ have the same vertices. 
\end{corollary}

\begin{proof}
This follows from Lemma~\ref{lem functor properties}(d) and Proposition~\ref{prop vertices}. The first statement also follows from \cite[Chapter~III, Lemma~4.14]{Feit1982}.
\end{proof}

Feit attributes the following theorem to Brauer.

\begin{theorem}\cite[Theorem 19.3]{Feit1982}
Let $F'$ be a finite field, $A$ a finite dimensional $F'$-algebra, $K/F'$ a field extension, $V$ an absolutely irreducible $K\otimes_{F'} A$-module such that $\mathrm{tr}_V(a)\in F'$ for every $a\in A$. Then $V$ has an $A$-form, i.e., there exists an absolutely irreducible $A$-module $W$ such that $K\otimes_{F'} W\cong V$ as $K\otimes_{F'} A$-modules.
\end{theorem}

\begin{corollary}\label{cor descent for irreducibles}
Let $V$ be an irreducible $FG$-module and let $F':=F^\Delta$ denote the fixed field of $\Delta:=\stab_\Gamma(V)$. Then there exists an (absolutely) irreducible $F'G$-module $W$ with $V\cong F\otimes_{F'} W$. 
\end{corollary}

\begin{proof}
This follows from the above theorem noting that if $\lexp{\sigma}{V}\cong V$ for some $\sigma\in\Gamma$ then $\sigma(\tr_V(g))=\tr_V(g)$ for all $g\in G$.
\end{proof}

Recall from \cite[Section 5.11]{Linckelmann2018a} that, for any field $F'$ of characteristic $p>0$, a {\em $p$-permutation $F'G$-module} $M$ is a direct summand of a finitely generated permutation $F'G$-module. Equivalently, the restriction of $M$ to any $p$-subgroup of $G$ is a permutation module. Also equivalently, the sources of the indecomposable direct summands of $M$ are trivial modules. We denote the Grothendieck group of $p$-permutation $F'G$-modules $V$ with respect to split short exact sequences by $T_{F'}(G)$. It is a commutative ring with multiplication induced by $-\otimes_{F'}-$. The class of $V$ in $T_{F'}(G)$ is denoted by $[V]$. The classes of indecomposable $p$-permutation $F'G$-modules form a standard $\ZZ$-basis of $T_{F'}(G)$.  If $E$ is a field extension of $F$ then the ring homomorphism $T_F(G)\to T_{E}(G)$ of scalar extension is an isomorphism, see \cite[Theorem~1.9]{BoltjeGlesser2007}. The Galois conjugation functors in (\ref{eqn Galois functors}) induce an action of the group $\Gamma$ on $T_F(G)$ via ring isomorphisms which stabilizes the standard basis. For a subfield $F'$ of $F$ the functors of scalar restriction and extension induce a group homomorphism $T_F(G)\to T_{F'}(G)$ and a ring homomorphism 
\begin{equation}\label{eqn T embedding}
  T_{F'}(G)\to T_F(G)
\end{equation}
which is injective by the Deuring-Noether Theorem and whose image is contained in the subring $T_F(G)^\Delta$ of $\Delta$-fixed points, where $\Delta:=\Gal(F/F')$. Note that also the abelian group $T_F(G)^\Delta$ has a standard $\ZZ$-basis, namely the $\Delta$-orbit sums of the standard basis of $T_F(G)$.
The goal of this section is the following theorem.

\begin{theorem}\label{thm Galois descent}
Let $F'$ be a subfield of $F$ and $\Delta:=\Gal(F/F')$. Then the ring homomorphism $T_{F'}(G)\to T_F(G)^\Delta$ in (\ref{eqn T embedding}) induced by scalar extension is an isomorphism, mapping the standard basis to the standard basis.
\end{theorem}

Before proving the above theorem we need the following proposition.

\begin{proposition}\label{prop descend to stabilizer}
Let $V$ be an indecomposable $p$-permutation $FG$-module and let $F'=F^{\Delta}$ be the fixed field of $\Delta:=\stab_{\Gamma}(V)$. Then there exists a unique (up to isomorphism) indecomposable $p$-permutation $F'G$-module $W$ such that $V\cong F\otimes_{F'} W$ as $FG$-modules. 
\end{proposition}

\begin{proof}
Let $P$ be a vertex of the module $V$. Then the Brauer construction $V(P)$ of $V$ is projective indecomposable as an $F[N_G(P)/P]$-module and its inflation is the Green correspondent of $V$, see \cite[Theorem~5.10.5]{Linckelmann2018a}. The quotient module $S:=V\left(P\right)/J\left(V\left(P\right)\right)$ is absolutely irreducible since $F$ is a splitting field.  Since the Green correspondence and taking projective covers commutes with Galois conjugation, the stabilizer of the isomorphism class of $S$ in $\Gamma$ is equal to $\Delta$.

By Corollary \ref{cor descent for irreducibles}, there exists an irreducible $F'[N_G(P)/P]$-module $T$ such that $S\cong F\otimes_{F'} T$. Let $W'$ be a projective indecomposable $F'[N_G(P)/P]$-module such that $W'/J(W')\cong T$, and let $W\in\lmod{F'G}$ be the Green correspondent of the inflation of $W'$. We will show that $V\cong F\otimes_{F'} W$.

First we claim that the projective $F[N_G(P)/P]$-module $F\otimes _{F'}W'$ is a projective cover of $S$. In fact,  $J(F\otimes_{F'}F'[N_G(P)/P]) = F\otimes_{F'} J(F'[N_G(P)/P])$, see \cite[Propositions 1.16.14 and 1.16.18]{Linckelmann2018a}, so that $J(F\otimes_{F'} W') = F\otimes_{F'} J(W')$. With this we obtain
\begin{equation*}
  (F\otimes_{F'} W')/J(F\otimes_{F'}W') = (F\otimes_{F'}W')/(F\otimes_{F'}J(W')) \cong F\otimes_{F'} (W'/J(W')) \cong S\,,
\end{equation*}
establishing the claim. Thus, $F\otimes_{F'} W'\cong V(P)$ as $F[N_G(P)/P]$-modules and also as $F[N_G(P)]$-modules after inflation. Since $W$ is the Green correspondent of $W'$, we have
\begin{equation*}
  F\otimes_{F'}W \mid F\otimes_{F'} \Ind_{N_G(P)}^G (W') \cong \Ind_{N_G(P)}^G(F\otimes_{F'} W') 
  \cong \Ind_{N_G(P)}^G(V(P))\,.
\end{equation*} 
But since the modules $V$ and $V(P)$ are Green correspondents, the module $V$ is the unique indecomposable direct summand of $\Ind_{N_G(P)}^G(V(P))$ with vertex $P$ and has multiplicity one in $\Ind_{N_G(P)}^G(V(P))$. Now Corollary~\ref{cor vertex} implies $F\otimes_{F'}W\cong V$, as desired.
\end{proof}

\bigskip\noindent
\begin{proof}{\em of Theorem~\ref{thm Galois descent}.}\quad It suffices to show that every standard basis element of $T_F(G)^\Delta$ comes via scalar extension from $T_{F^{\Delta}}(G)$. So let $V$ be an indecomposable $p$-permutation $FG$-module and set $\Delta':=\stab_\Gamma([V])$. Then the $\Delta$-orbit sum of $[V]$, i.e., the class of $\bigoplus_{\sigma\in\Delta/(\Delta\cap \Delta')} \lexp{\sigma}{V}$ is a standard basis element of $T_F(G)^\Delta$ and every standard basis element is of this form. By Proposition~\ref{prop descend to stabilizer} there exists an indecomposable $p$-permutation $F^{\Delta'}G$-module $W'$ such that $V\cong F\otimes_{F^{\Delta'}} W'$. By Lemma~\ref{lem functor properties}(d), we have 
\begin{equation*}
   \bigoplus_{\sigma\in\Delta/(\Delta\cap \Delta')} \lexp{\sigma}{V}
   \cong F\otimes_{F^{\Delta'}}\bigl(\bigoplus_{\sigma\in\Delta\Delta'/\Delta'} \lexp{\sigma}{W'}\bigr) 
   \cong F\otimes_{F^{\Delta\Delta'}} (W'_{F^{\Delta\Delta'}}) \cong F\otimes_{F^\Delta} W
\end{equation*}
with $W:= F^{\Delta}\otimes_{F^{\Delta\Delta'}} (W'_{F^{\Delta\Delta'}})$.
\end{proof}

%%%%%%%%%%%%%%%%%% SECTION 3 %%%%%%%%%%%%%%%%%%%%%%%%%%%%%%%%%

\section{Endosplit $p$-permutation resolutions}\label{sec resolutions}

In this section, and only this section, $F$ can be any field of characteristic $p$. The following concept is due to Rickard, see \cite[Section 7]{Rickard1996}.

\begin{definition}
Let $M$ be a finitely generated $FG$-module. An {\em endosplit $p$-permutation resolution} of $M$ is a bounded chain complex $X$ of $p$-permutation $FG$-modules with homology concentrated in degree $0$  such that $H_0(X)\cong M$ and such that $X\otimes_F X^\circ$ is split as chain complex of $FG$-modules (with $G$ acting diagonally and $X^\circ$ denoting the $F$-dual of $X$). Here, $X^\circ$ is again considered as a chain complex.
\end{definition}

\begin{Remark}\label{rem endosplit}
Let $X$ be an endosplit $p$-permutation resolution of a finitely generated $FG$-module $M$. 

\smallskip
(a) Every direct summand $X'$ of $X$ is again an endosplit $p$-permutation resolution of $H_0(X')$. 

\smallskip
(b) We can decompose $X$ into a direct sum $X=X'\oplus X''$ of chain complexes such that $X''$ is contractible and $X'$ has no contractible non-zero direct summand. With $X$, also $X'$ is then an endo-split $p$-permutation resolution of $M$. If $X''=0$, we say that $X$ is {\em contractible-free}. 

\smallskip
(c) Taking the $0$-th homology induces an $F$-algebra isomorphism
\begin{align}\label{isom markus}
\rho: \End_{K(\lmod{FG})}(X)\cong \End_{FG}(M)\,,
\end{align} 
where $K(\lmod{FG})$ denotes the homotopy category of chain complexes in $\lmod{FG}$, see \cite[Proposition 7.11.2]{Linckelmann2018b}. If $N\mid M$, then the projection map onto $N$ yields an idempotent in $\End_{FG}(M)$ and hence an idempotent in $\End_{K(\lmod{FG})}(X)$ via the isomorphism in (\ref{isom markus}). This idempotent lifts to an idempotent $\alpha$ in $\End_{Ch(\lmod{FG})}(X)$, where $Ch(\lmod{FG})$ denotes the category of chain complexes in $\lmod{FG}$. It follows that the direct summand $\alpha(X)$ of $X$ is an endosplit $p$-permutation resolution of $N$. The lifted idempotent is not unique up to conjugation, but $\alpha(X)$ is unique up to isomorphism and contractible direct summands. Therefore, if $X$ is contractible-free, then $\alpha(X)$ is uniquely determined by $N$ up to isomorphism in $Ch(\lmod{FG})$.

\smallskip
(d) Suppose that $M\cong F\otimes_{F'} M'$ for some subfield $F'\subseteq F$ and some $M'\in\lmod{F'G}$. Then, $M_{F'}\cong M'^{[F:F']}$ in $\lmod{F'G}$ and $X_{F'}\in Ch(\lmod{F'G})$ is an endosplit $p$-permutation resolution of $M'^{[F:F']}$. By Part~(c), also $M'$ has an endosplit $p$-permutation resolution. Conversely, if $M'$ has an endosplit $p$-permutation resolution $X'\in Ch(\lmod{F'G})$ then $F\otimes_{F'}X'$ is an endosplit $p$-permutation resolution of $F\otimes_{F'}M'\cong M$.
\end{Remark}

\begin{lemma}\label{diagram}
Let $X_V$, $X_U$, $X_{V'}$ and $X_{U'}$ be endosplit $p$-permutation resolutions of  $V,U,V',U'\in\lmod{FG}$, respectively, and assume that $X_V$ and $X_{V'}$ are contractible-free. Suppose further that $X_V\oplus X_U\cong X_{V'} \oplus X_{U'}$ in $Ch(\lmod{FG})$ are endo-split $p$-permutation resolutions of $V\oplus U$ and $V'\oplus U'$, respectively, and that $V\cong V'$ in $\lmod{FG}$. Then $U\cong U'$ in $\lmod{FG}$ and $X_V\cong X_{V'}$ in $Ch(\lmod{FG})$.
%
%Then by Lemma \ref{lem directsummandresolution} $X_W=X_V\oplus X_U$ for some endosplit $p$-permutation resolutions $X_U$ and $X_V$ of $U$ and $V$, respectively. Similarly,  $X_{W'}=X_{V'}\oplus X_{U'}$ for some endosplit $p$-permutation resolutions $X_{U'}$ and $X_{V'}$ of $U'$ and $V'$, respectively. Assume that we have $V\cong V'$ as $FG$-modules and $X_W\cong X_{W'}$ as complexes of $FG$-modules. Then we also have $W\cong W'$ as $FG$-modules and $X_V\cong X_{V'}$ as complexes of $FG$-modules.
\end{lemma}

\begin{proof}
Taking $0$-th homology of $X_V\oplus X_U$ and $X_{V'}\oplus X_{U'}$ yields $V\oplus U\cong V'\oplus U'$, and the Krull-Schmidt Theorem implies $U\cong U'$.  For the second statement let $\phi\colon X_V\oplus X_U \to X_{V'}\oplus X_{U'}$ be an isomorphism in $Ch(\lmod{FG})$. Then $\phi(X_V)$ and $X_{V'}$ are both direct summands of $X_{V'}\oplus X_{U'}$ and contractible-free endo-split $p$-permutation resolutions of $V$. Therefore, by \cite[Proposition~7.11.2]{Linckelmann2018b} (see also Remark~\ref{rem endosplit}(b)) they are isomorphic.
\end{proof}

%%%%%%%%%%%%%%%%%%%%% SECTION 4 %%%%%%%%%%%%%%%%%%%%%%%%%%%%%%%%%%%

\section{Galois descent of $p$-permutation equivalences}\label{sec descent of equivalences}

Since the Galois group $\Gamma$ acts via $\FF_p$-algebra automorphisms on the group algebra $FG$ and also on $Z(FG)$, it permutes the block idempotents of $FG$.

\begin{proposition}\cite[Proposition 4.1]{BKY2020}\label{prop Gamma and blocks}
{\rm (a)} Let $b$ be a block idempotent of $FG$. Then $\btilde:=\Tr_{\Gamma}(b):=\sum_{\sigma\in\Gamma/\stab_{\Gamma}(b)} \lexp{\sigma}{b}$ is a block idempotent of $\FF_pG$.

\smallskip
{\rm (b)} The map $b\mapsto \btilde$ induces a bijection between the set of $\Gamma$-orbits of block idempotents of $FG$ and the set of block idempotents of $\FF_pG$.

\smallskip
{\rm (c)} If $b$ is a block idempotent of $FG$ and $\btilde:=\Tr_\Gamma(b)$ is the block idempotent of $\FF_pG$ associated to it as in (a) then $\btilde$ and $b$ have the same defect groups.
\end{proposition}

\begin{lemma}\label{lem brauercorrespondentblocks}
Let $b$ be a block of $FG$ with a defect group $P$ and $c$ be the block of $FN_G(P)$ which is in Brauer correspondence with $b$. For any $\sigma\in \Gamma$, the blocks $\lexp{\sigma}{b}$ and $\lexp{\sigma}{c}$ are again in Brauer correspondence. In particular, the stabilizers of $b$ and $c$ in $\Gamma$ are the same. Moreover, the blocks $\btilde=\Tr_{\Gamma}(b)$ of $\FF_pG$ and $\ctilde=\Tr_{\Gamma}(c)$ of $\FF_pN_G(P)$ are Brauer correspondents. 
\end{lemma}

\begin{proof}
The first assertion follows immediately from the fact that the action of $\Gamma$ and the Brauer map $\Br_P$ commute. We have $\sigma\in \mathrm{stab}_{\Gamma}(b)$ $\iff$ $\sigma(b)=b$ $\iff$ $\Br_P(\sigma(b))=\Br_P(b)$ $\iff$ $\sigma(c)=c$ $\iff$ $\sigma\in\mathrm{stab}_{\Gamma}(c)$. The last statement follows easily from the additivity of the Brauer map. 
\end{proof}

Let $F'$ be a field of characteristic $p$ and let $b$ and $c$ be central idempotents of $F'G$ and $F'H$, respectively. As usual we identify $F'[G\times H]=F'G\otimes_{F'} F'H$ as $F$-algebras and we identify $(F'Gb,F'Hc)$-bimodules with left $F'[G\times H](b\otimes c^*)$-modules, where $c^*$ is defined by applying the $F'$-linear extension of $h\mapsto h^{-1}$ to $c$. We write $T^\Delta(F'Gb,F'Hc)$ for the subgroup of $T_{F'}(G\times H)$ spanned by indecomposable $F'[G\times H](b\otimes c^*)$-modules whose vertices are twisted diagonal, i.e., of the form $\{(\phi(y),y)\mid y\in Q\}$ for some isomorphism $\phi\colon Q\to P$ between $p$-subgroups $P$ and $Q$ of $G$ and $H$, respectively. Recall from \cite{BoltjePerepelitsky2020} that a {\em $p$-permutation equivalence} between $F'Gb$ and $F'Hc$ is an element $\omega\in T^\Delta(F'Gb, F'Hc)$ such that $\omega\cdot_{H} \omega^\circ = [F'Gb]$ in $T^\Delta(F'Gb, F'Gb)$ and $\omega^\circ\cdot_G\omega=[F' Hc]$ in $T^\Delta(F'Hc,F'Hc)$. Here, $\cdot_H$ is induced by $-\otimes_{F'H}-$, and $\omega^\circ\in T_{F'}(H\times G)$ is given by taking the $F'$-dual of $\omega$. Note that if $F'=F$ then $\stab_{\Gamma}(\omega)\le\stab_{\Gamma}(b)$ and $\stab_{\Gamma}(\omega)\le\stab_{\Gamma}(c)$.

\begin{lemma}\label{lem descending}
Let $b$ and $c$ be block idempotents of $FG$ and $FH$, respectively. Let $\btilde$ and $\ctilde$ denote the block idempotents of $\FF_pG$ and $\FF_pH$ associated to $b$ and $c$, respectively, as in Proposition~\ref{prop Gamma and blocks}(a). Moreover, let $\omega\in T^\Delta(FGb,FHc)$ be a $p$-permutation equivalence between $FGb$ and $FHc$. Suppose that we have $\Delta:=\stab_{\Gamma}(\omega)=\stab_{\Gamma}(b)=\stab_{\Gamma}(c)$. Then there exists a $p$-permutation equivalence between $\FF_pG\btilde$ and $\FF_pH\ctilde$.
\end{lemma}

\begin{proof}
For any $\sigma\in\Gamma$, the Galois conjugate $\lexp{\sigma}\omega$ is a $p$-permutation equivalence between $FG\lexp{\sigma}b$ and $FH\lexp{\sigma}c$. Hence the sum $\sum_{\sigma\in\Gamma/\Delta}\lexp{\sigma}\omega\in T^\Delta(FG\btilde,FH\ctilde)$ is a $p$-permutation equivalence between $\bigoplus_{\sigma\in\Gamma/\Delta}FG\lexp{\sigma}b=FG\btilde$ and $\bigoplus_{\sigma\in\Gamma/\Delta}FH\lexp{\sigma}c=FH\ctilde$. Note that the sum $\sum_{\sigma\in\Gamma/\Delta}\lexp{\sigma}\omega\in T_F(G\times H)$ is fixed under $\Gamma$. By Theorem~\ref{thm Galois descent}, there exists $\tilde{\omega}\in T_{\FF_p}(G,H)$ such that $\sum_{\sigma\in\Gamma/\Delta}\lexp{\sigma}\omega=F\otimes_{\FF_p}\tilde{\omega}$. It follows that $\tilde{\omega}$ is a $p$-permutation equivalence between $\FF_pG\btilde$ and $\FF_pH\ctilde$. 
\end{proof}

%%%%%%%%%%%%%%%%%% SECTION 5 %%%%%%%%%%%%%%%%%%%%%%%%%%%%%

\section{$p$-nilpotent groups}\label{sec $p$-nilpotent groups}

Throughout this section we assume that $G$ is a $p$-nilpotent group. Thus, $G$ has a normal $p'$-subgroup $N$ such that $G/N$ is a $p$-group. We fix a block idempotent $b$ of $FG$ and denote by $\btilde:=\Tr_\Gamma(b)$ the corresponding block idempotent of $\FF_p G$, see Proposition~\ref{prop Gamma and blocks}(a). 
Moreover, we fix a block idempotent $e$ of $FN$ such that $be\neq 0$. Then $b=\sum_{g\in G/S} \lexp{g}{e}$, where $S:=\stab_G(e)$, and the idempotent $e$ is also a block idempotent of $FS$. 
Let $Q$ be a Sylow $p$-subgroup of $S$. Then $Q$ is a defect group of the block idempotents $e$ of $FS$, $b$ of $FG$, and $\btilde$ of $\FF_pG$.
Finally, set $\etilde:=\Tr_\Gamma(e)$, the block idempotent of $\FF_pN$ determined by $e$ and set $\Stilde:=\stab_G(\etilde)$. Then $S\le \Stilde$ and $\btilde=\sum_{G/\Stilde} \lexp{g}{\etilde}$.

\smallskip
The group $\Gamma\times G$ acts on the block idempotents of $FN$. Set $X:=\stab_{\Gamma\times G}(e)$. Since $\stab_G(e)=S$ we have $k_2(X):=\{g\in G\mid (1,g)\in X\}=S$. Similarly, $k_1(X):=\{\sigma\in\Gamma\mid (\sigma,1)\in X\}=\stab_\Gamma(e)$. Next we determine the images of $X$ under the projection maps $p_1\colon \Gamma\times G\to\Gamma$ and $p_2\colon \Gamma\times G\to G$.

\begin{lemma}\label{lem stabilizers 1}
One has $p_2(X)=\Stilde$ and $S\trianglelefteq\Stilde$.
\end{lemma}
\begin{proof}
Let $g\in p_2(X)$. There exists $\sigma\in\Gamma$ such that $\lexp{(\sigma,g)}e=e$. Therefore we have
\begin{equation*}
  \etilde=\Tr_{\Gamma}(e)=\Tr_{\Gamma}(\lexp{(\sigma, g)}e)=\Tr_{\Gamma}(\lexp{\sigma}(\lexp{g}e))=
  \Tr_{\Gamma}(\lexp{g}e)=\lexp{g}\Tr_{\Gamma}(e)=\lexp{g}\etilde\,.
\end{equation*}
This shows that $g\in\stab_G(\etilde)=\Stilde$ and hence that $p_2(X)\le \Stilde$. 

Now let $\stilde\in \Stilde$. Then
\begin{equation*}
\Tr_{\Gamma}(\lexp{\stilde}e)=\lexp{\stilde}(\Tr_{\Gamma}(e))=\lexp{\stilde}\etilde=\etilde\,.
\end{equation*}
Since the blocks $e$ and $\lexp{\stilde}e$ have the same Galois trace, they must be $\Gamma$-conjugate, and therefore $\stilde\in p_2(X)$. This proves the first statement. The second statement holds, since $k_2(X)$ is normal in $p_2(X)$ in general, see \cite[p.~24]{Bouc2010}.
\end{proof}

Next, set $e':=\sum_{\stilde\in\Stilde/S}\lexp{\stilde}{e}$. Then $e'$ is a block idempotent of $F\Stilde$ and $b=\sum_{g\in G/\Stilde}\lexp{g}{e'}$.

\begin{lemma}\label{lem stabilizers 2}
One has $\stab_{\Gamma}(e')=\stab_{\Gamma}(b)=p_1(X)$. Moreover, $\Stilde/S\cong \stab_\Gamma(b)/\stab_\Gamma(e)$ is cyclic.
\end{lemma}

\begin{proof}
We have $\stab_{\Gamma}(e')\le\stab_{\Gamma}(b)$ since $b=\sum_{g\in G/\Stilde}\lexp{g}{e'}$. 
Next, let $\sigma\in p_1(X)$. Then there exists $\stilde_0\in \Stilde$ such that $(\sigma,\stilde_0)\in X$, and
\begin{equation*}
  e'=\sum_{\stilde\in \Stilde/S}\lexp{\stilde}e=\sum_{\stilde\in \Stilde/S}\lexp{\stilde}(\lexp{(\sigma,\stilde_0)}e)
  =\lexp{\sigma}(\sum_{\stilde\in \Stilde/S}\lexp{\stilde}(\lexp{\stilde_0}e))=\lexp{\sigma}e'\,,
\end{equation*}
where the last equation holds, because $S\trianglelefteq \Stilde$. This shows that $\sigma\in \stab_{\Gamma}(e')$ and hence $p_1(X)\le\stab_{\Gamma}(e')$. 
Finally, let $\sigma\in\stab_{\Gamma}(b)$. Then $\lexp{\sigma}{b}=b$ implies that $ \lexp{\sigma}(\sum_{g\in G/S}\lexp{g}e) = \sum_{g\in G/S}\lexp{g}e$.
Therefore there exists $g\in G$ such that $\lexp{(\sigma, g)}e=e$, i.e., $\sigma\in p_1(X)$. The proof of the first statement is now complete. The second statement follows from the general isomorphism $p_1(X)/k_1(X)\cong p_2(X)/k_2(X)$, see \cite[p.~24]{Bouc2010}, and since $\Gamma$ is cyclic.
\end{proof}

\begin{lemma}\label{lem e' and etilde}
One has $\Tr_{\Gamma}(e')=\etilde$. In particular, $\etilde$ is a block idempotent of $\FF_p\Stilde$. 
\end{lemma}
\begin{proof}
By Lemma \ref{lem stabilizers 2} and since $k_1(X)=\stab_\Gamma(e)$, we have
\begin{align*}
  \Tr_{\Gamma}(e')&=\sum_{\sigma\in\Gamma/p_1(X)}\lexp{\sigma}e'=
       \sum_{\sigma\in\Gamma/p_1(X)}\sum_{\stilde\in \Stilde/S}\lexp{\sigma}(\lexp{\stilde}e)
       =\sum_{\sigma\in\Gamma/p_1(X)}\sum_{\tau\in p_1(X)/k_1(X)}\lexp{\sigma}(\lexp{\tau}e)\\
       & =\sum_{\sigma\in\Gamma/k_1(X)}\lexp{\sigma}e=\Tr_{\Gamma}(e)=\etilde\,,
\end{align*}
as desired. The third equation holds, since the classes of $\stilde$ and $\tau$ correspond under the isomorphism $p_2(X)/k_2(X)\cong p_1(X)/k_1(X)$ if and only if $(\tau,\stilde)\in X$, see \cite[p.~24]{Bouc2010}.
\end{proof}

Let $V$ denote the unique (up to isomorphism) simple $FNe$-module. By Theorem~\ref{thm Galois descent} and since $\stab_\Gamma(V)=\stab_\Gamma(e) = k_1(X)$, there exists a unique simple $\FF_pN$-module $\Vtilde$ such that 
\begin{equation}\label{eqn Vtilde}
   \bigoplus_{\sigma\in \Gamma/k_1(X)} \lexp{\sigma}{V}\cong F\otimes_{\FF_p}\Vtilde\,.
\end{equation} 
Since $\etilde$ acts as identity on the above direct sum, $\Vtilde$ is a simple $\FF_p N\etilde$-module. Since $V$ is absolutely irreducible, it extends to a (unique up to isomorphism) simple $FSe$-module which we denote again by $V$. Similarly, each $\lexp{\sigma}{V}$ can be viewed as $FS$-module, so that the left hand side in (\ref{eqn Vtilde}) has an $F S\etilde$-module structure and is $\Gamma$-invariant. Again, by Theorem~\ref{thm Galois descent}, the left hand side in (\ref{eqn Vtilde}) regarded as $FS$-module has an $\FF_p$-form $W\in\lmod{\FF_p S\etilde}$. Restriction from $S$ to $N$ and the Deuring-Noether Theorem then imply that $\Res^S_N(W) \cong \Vtilde$. Thus, $\Vtilde$ extends to a simple $FS\etilde$-module and (\ref{eqn Vtilde}) is an isomorphism of $FS\etilde$-modules.

\begin{proposition}\label{prop extendingsimple}
The $\FF_pS\etilde$-module $\Vtilde$ extends to an $\FF_p\Stilde\etilde$-module.
\end{proposition}

\begin{proof}
By Fong's first reduction theorem, $W:=\Ind_S^{\Stilde} V$ is the unique simple $F\Stilde e'$-module (up to isomorphism) and $\stab_{\Gamma}(W)=\stab_{\Gamma}(e')=p_1(X)$. By Theorem~\ref{thm Galois descent}, there exists a simple $\FF_p\Stilde$-module $\Wtilde$ such that $\bigoplus_{\sigma\in \Gamma/p_1(X)}\lexp{\sigma}{W}\cong F\otimes_{\FF_p}\Wtilde$. Restriction to $S$ implies
\begin{align*}
   \Res^{\Stilde}_S(F\otimes_{\FF_p}\Wtilde) & \cong \Res^{\Stilde}_S(\bigoplus_{\sigma\in \Gamma/p_1(X)}\lexp{\sigma}{W})
         \cong \bigoplus_{\sigma\in \Gamma/p_1(X)}\lexp{\sigma}{(\Res^{\Stilde}_S W)}
         \cong \bigoplus_{\sigma\in \Gamma/p_1(X)}\lexp{\sigma}{(\Res^{\Stilde}_S\Ind^{\Stilde}_S V)}\\
         & \cong \bigoplus_{\sigma\in \Gamma/p_1(X)}\lexp{\sigma}{(\bigoplus_{\stilde\in \Stilde/S}\lexp{\stilde}{V})}
         \cong \bigoplus_{\sigma\in \Gamma/k_1(X)}\lexp{\sigma}{V} \cong F\otimes_{\FF_p}\Vtilde\,,
\end{align*}
since $\bigoplus_{\stilde\in \Stilde/S}\lexp{\stilde}{V}\cong \bigoplus_{\tau\in p_1(X)/k_1(X)}\lexp{\tau}{V}$, which follows from the argument at the end of the proof of the previous proposition. This shows that $\Res^{\Stilde}_S\Wtilde=\Vtilde$ and the result follows. 
\end{proof}

\begin{Remark}
Proposition \ref{prop extendingsimple} extends the results of Michler \cite[Theorem 3.7]{Michler1973} (z=1 in part(e)). 
\end{Remark}

Now set $H:=N_G(Q)$, which is again a $p$-nilpotent group, and set $M:=O_{p'}(H)$, the largest normal $p'$-subgroup of $H$. Then 
\begin{align*}
M=H\cap N= C_N(Q)\,.
\end{align*} 
Let $c$ denote the block idempotent of $FH$ which is in Brauer correspondence with $b$. Then, by Lemmas~\ref{lem brauercorrespondentblocks} and \ref{lem stabilizers 2}, $\stab_{\Gamma}(c)=\stab_{\Gamma}(b)=p_1(X)$ and $\ctilde:=\Tr_\Gamma(c)$ is the Brauer correspondent of $\btilde$.

Further, let $f$ denote the block idempotent of $FM$ whose irreducible module is the Glaubermann correspondent of the $Q$-stable irreducible module $V\in\lmod{FN}$. Then $f$ is $QM=QC_N(Q)=N_S(Q)=N_G(Q)\cap S=H\cap S=:T$-stable and hence it remains a block idempotent of $FT$. By \cite{Alperin1976}, the block idempotents $e$ of $FS$ and $f$ of $FT$ are Brauer correspondents.

\begin{lemma}
One has $c=\Tr_T^H(f)$ and $\stab_H(f)=T$.
\end{lemma}
\begin{proof}
Since the block idempotents $b$ and $c$ are Brauer correspondents, we have
\begin{align*}
   c & =\Br_Q(b) = \Br_Q(\Tr_{S}^G(e)) = 
      \Br_Q\big(\sum_{x\in Q\setminus G/S}\Tr^Q_{Q\cap \lexp{x}S}(\lexp{x}e)\big)\\
   &  = \sum_{x\in Q\setminus G/S}\Br_Q\big(\Tr^Q_{Q\cap \lexp{x}S}(\lexp{x}e)\big)
        = \sum_{\substack{x\in Q\setminus G\\ Q\le\lexp{x}S}}\Br_Q(\lexp{x}e)\,.
\end{align*}
The condition $Q\le\lexp{x}S$ implies that $\lexp{x^{-1}}Q\le S$ and hence $\lexp{x^{-1}}Q=\lexp{s}Q$ for some $s\in S$ since $Q$ is a Sylow $p$-subgroup of $S$. This means that $xs\in N_G(Q)$ and so $x\in N_G(Q)S$. Therefore the above sum can be written as
\begin{align*}
c=\sum_{x\in N_G(Q)/(N_G(Q)\cap S)}\Br_Q(\lexp{x}e)=\sum_{x\in H/T}\lexp{x}{\Br_Q(e)}=\Tr_T^H(f)\,,
\end{align*}
since $f=\Br_Q(e)$. This proves the first assertion. The group $\stab_H(f)$ has the group $Q$ as a Sylow $p$-subgroup, since $Q$ is a defect group of the bock $c$. This shows that $\stab_H(f)=QM=T$, as desired. 
\end{proof}

Let $\ftilde:=\Tr_{\Gamma}(f)$, $\Ttilde:=\stab_H(\ftilde)$, $f':=\Tr_T^{\Ttilde}(f)$ and $Y:=\stab_{\Gamma\times H}(f)$. Since the blocks $e$ and $f$ are Brauer correspondents, we have 
\begin{equation}\label{eqn k_1(X)=k_1(Y)}
   k_1(Y)=\stab_{\Gamma}(f)=\stab_{\Gamma}(e)=k_1(X)\,.
\end{equation} 
Moreover, by Lemmas~\ref{lem stabilizers 1} and \ref{lem stabilizers 2}, 
\begin{gather}\label{eqn p_1(X)=p_1(Y)}
   \stab_\Gamma(f') = p_1(Y) = \stab_{\Gamma}(c) = \stab_\Gamma(b) = p_1(X) = \stab_\Gamma(e') \,,\\
   p_2(Y)=\Ttilde \quad  \text{and} \quad k_2(Y)=T\,,\notag
\end{gather}
and therefore
\begin{equation}\label{eqn S versus T}
   \Ttilde/T = p_2(Y)/k_2(Y) \cong p_1(Y)/k_1(Y) = p_1(X)/k_1(X) \cong p_2(X)/k_2(X) = \Stilde/S 
\end{equation}
which implies that $\Ttilde=H\cap \Stilde$. 

\begin{nothing}\label{noth Rickard}
We recall Rickard's construction of a {\em splendid Rickard equivalence} between $FSe$ and $FQ$, i.e., a bounded chain complex $X$ of relatively $\Delta(Q)$-projective $p$-permutation $(FSe,FQ)$-bimodules such that $X\otimes_{FQ}X^\circ\cong FSe$ and $X^\circ\otimes_{FS}X\cong FQ$ in the homotopy categories of $(FSe,FSe)$-bimodules and $(FQ,FQ)$-bimodules, respectively, where $FSe$ and $FQ$ are considered as chain complexes concentrated in degree $0$. For more details we refer the reader to \cite{Rickard1996}. 

Set $\Delta_QS:=\{(nq,q): n\in N, q\in Q\}\le S\times Q$ and note that $p_1\colon S\times Q\to S$ restricts to an isomorphism $\Delta_QS\myiso S$. 
The module $\Res^S_Q V$ is a capped endopermutation $FQ$-module. In everything that follows, we suppose that
\begin{equation}\label{eqn hypothesis}
   \text{{\em $\Res^S_Q V$ has an endosplit $p$-permutation resolution $X_V$.}}
\end{equation}
By the proof in \cite[Lemma~7.7]{Rickard1996}, see also Remark~(a) at the end of Section~7 in \cite{Rickard1996}, the induced complex $\Ind^S_Q(X_V)$ is an endosplit $p$-permutation resolution of $\Ind^S_Q\Res^S_Q V$ as $FS$-modules. Since $V \mid \Ind^S_Q\Res^S_Q V$, there exists a direct summand $Y_{V}$ of $\Ind^S_Q X_V$ such that $Y_{V}$ is an endosplit $p$-permutation resolution of $V$ as an $FS$-module, and we may choose $Y_V$ to be contractible-free, see Remark~\ref{rem endosplit}(c) and (b). The induced chain complex $\Ind_{\Delta_QS}^{S\times Q} Y_{V}$ is then a splendid Rickard equivalence between $FSe$ and $FQ$, see \cite[Theorem~7.8]{Rickard1996} and its subsequent Remark~(a).
\end{nothing}

\begin{proposition}\label{prop Rickard equivalence over Fp}
Suppose that $\stab_{\Gamma}(e)=\stab_{\Gamma}(b)$ and that $\Res^S_QV$ has an endosplit $p$-permutation resolution.

\smallskip
{\rm (a)} There exists a splendid Rickard equivalence between $\FF_pS\etilde$ and $\FF_pT\ftilde$.

\smallskip
{\rm (b)} There exists a splendid Rickard equivalence between $\FF_pG\btilde$ and $\FF_pH\ctilde$. 
\end{proposition}

\begin{proof}
(a) The equality $\stab_{\Gamma}(e)=\stab_{\Gamma}(b)$ implies that we have 
\begin{equation*}
\stab_{\Gamma}(f)=\stab_{\Gamma}(c)=\stab_{\Gamma}(b)=\stab_{\Gamma}(e) \,, \quad \Stilde=S\,, \quad  \Ttilde=T \,, \quad e'=e \quad \text{and} \quad f'=f\,,
\end{equation*} 
by (\ref{eqn k_1(X)=k_1(Y)}) and (\ref{eqn p_1(X)=p_1(Y)}). Let $\FF_p[e]$ denote the smallest field containing the coefficients of the idempotent $e$. Then $F':=\FF_p[e]=\FF_p[f]\subseteq F$. By Corollary~\ref{cor descent for irreducibles}, there exists an absolutely simple $F'N$-module $V'$ such that $V\cong F\otimes_{F'}V'$, the unique simple module in the block $FNe$. Since $e$ is $S$-stable, also $V'$ extends to an $F'S$-module that we again denote by $V'$. Then $V\cong F\otimes_{F'}V'$ also as $FSe$-modules. Since $V\in\lmod{FS}$ has an endosplit $p$-permutation resolution, also $V'\in\lmod{F'S}$ has an endosplit $p$-permutation resolution $X'\in Ch(\lmod{F'S})$, see Remark~\ref{rem endosplit}(d). Since $F'$ is a splitting field of $V'$ as $F'N$-module, we may use the results from Theorem~7.8 and its subsequent Remark~(a) in \cite{Rickard1996} in order to see that $\Ind_{\Delta_QS}^{S\times Q}(X')$ is a splendid Rickard equivalence between $F'Se$ and $F'Q$. Using the unique $F'$-form $U'\in\lmod{F'M}$ of $U\in \lmod{FM}$ and its unique extension to an $F'T$-mdoule, we similarly obtain that $\Ind_{\Delta_QT}^{T\times Q}(U')$ induces a splendid Morita equvivalence between $F'Tf$ and $F'Q$. Thus, the chain complex
\begin{equation*}
Z'=\Ind_{\Delta_QS}^{S\times Q} (X')\otimes_{FQ} \Ind_{(\Delta_QT)^\circ}^{Q\times T} (U')^\circ
\end{equation*}
is a splendid Rickard equivalence between $F'Se$ and $F'Tf$. The result now follows from \cite[Theorem 6.5]{KL2018}.

\smallskip
(b) The $p$-permutation bimodule $\FF_pG\etilde$ induces a Morita equivalence, hence a splendid Rickard equivalence, between $\FF_pG\btilde$ and $\FF_p S\etilde$. Similarly, the bimodule $\FF_pH\ftilde$ induces a splendid Rickard equivalence between $\FF_p H\ctilde$ and $\FF_p T\ftilde$. The result follows now from Part~(a).
\end{proof}

\begin{Remark}\label{rem on hypo A}
By the classification of indecomposable capped endopermutation modules, the hypothesis in (\ref{eqn hypothesis}) is satisfied if $p$ is odd, or if $p=2$ and $Q$ does not have a subquotient isomorphic to the quaternion group of order $8$. See \cite{Thevenaz2007} for more details. Therefore, Theorem~A follows from Proposition~\ref{prop Rickard equivalence over Fp}.
\end{Remark}

\begin{nothing}\label{noth Z}
(a) By Fong's first reduction theorem, the $(F\Stilde e', FSe)$-bimodule $F\Stilde e$ induces a Morita equivalence between $F\Stilde e'$ and $FSe$. Hence the complex $F\Stilde e\otimes_{FS} \Ind_{\Delta_QS}^{S\times Q} Y_{V}$ gives a splendid Rickard equivalence between $F\Stilde e'$ and $FQ$.

\smallskip
(b) For any $FS$-module $M$, let $M\otimes_F FQ$ be the $(FS, FQ)$-bimodule, given by 
\begin{equation*}
   \text{$s (m\otimes x) y := sm\otimes qxy$, \quad for $s=nq\in S$, $n\in N$, $m\in M$, and $x,y,q\in Q$.}
\end{equation*}
It is straightforward to check that the map
\begin{align*}
  \phi_M: M\otimes_F FQ &\to  F[S\times Q]\otimes_{F\Delta_QS} M \\
  v\otimes q &\mapsto (1, q^{-1})\otimes v
\end{align*} 
is an isomorphism of $(FS, FQ)$-bimodules and that it is natural in $M$. 
Therefore, it yields an isomorphism
\begin{equation*}
   Y_V\otimes_F FQ \cong \Ind_{\Delta_QS}^{S\times Q} Y_{V}
\end{equation*}
of chain complexes of $(FSe, FQ)$-bimodules.

\smallskip
(c) Let $U$ be the simple $FM$-module belonging to the block idempotent $f$. Since $Q$ is normal in $H$, we have $T=Q\times M$. Thus, the unique extension of $U$ to $T$ (with $Q$ acting trivially on $U$) is a $p$-permutation $FT$-module and plays the same role as the complex $Y_V$. Similar as in (a), the bimodule $F\Ttilde f\otimes_{FT} \Ind_{\Delta_QT}^{T\times Q} U$ induces a splendid Rickard equivalence between $F\Ttilde f'$ and $FQ$.  

\smallskip
(d) Altogether, the complex 
\begin{equation*}
   Z:=F\Stilde e\otimes_{FS} \Ind_{\Delta_QS}^{S\times Q} Y_{V}\otimes_{FQ} \Ind_{(\Delta_QT)^\circ}^{Q\times T} (U)^\circ
   \otimes_{FT} fF\Ttilde
\end{equation*} 
induces a splendid Rickard equivalence between $F\Stilde e'$ and $F\Ttilde f'$. Here, $(\Delta_QT)^\circ:=\{(q,t)\in Q\times T\mid (t,q)\in\Delta_QT\}$. Set 
\begin{equation*}
   \omega:=\sum_{n\in \ZZ} (-1)^n[Z_n]\in T^\Delta(F\Stilde e', F\Ttilde f')\,. 
\end{equation*}
By \cite[Theorem~1.5]{BX2008}, $\omega$ is a $p$-permutation equivalence between $F\Stilde e'$ and $F\Ttilde f'$. Moreover, the isomorphism in (b) implies that
\begin{align*}
   Z&= F\Stilde e\otimes_{FS} \Ind_{\Delta_QS}^{S\times Q} Y_{V}\otimes_{FQ} 
         \Ind_{(\Delta_QT)^\circ}^{Q\times T} (U)^\circ\otimes_{FT} fF\Ttilde \\
   & \cong F\Stilde e\otimes_{FS} \left(Y_V\otimes_F FQ\right)\otimes_{FQ} \left(FQ\otimes_F U^\circ \right)\otimes_{FT} 
       fF\Ttilde \\
   & \cong F\Stilde e\otimes_{FS} \left(Y_V\otimes_F FQ\otimes_F U^\circ \right)\otimes_{FT} fF\Ttilde \,.
\end{align*}
\end{nothing}

Let $R$ be a Sylow $p$-subgroup of $\Ttilde$ containing $Q$. Then $\Ttilde=RM$, and, by (\ref{eqn S versus T}), $R$ is also a Sylow $p$-subgroup of $\Stilde$ so that $\Stilde = RN$.

\smallskip
The following diagram depicts the subgroups, block idempotents, and modules introduced so far.

\bigskip
\bigskip
\bigskip
\bigskip

\begin{center}
\unitlength 7mm
{\scriptsize
\begin{picture}(7,7)
\put(6,9){$\bullet$}  \put(6.2,9.2){$G, \btilde, b$}
\put(6,7){$\bullet$}  \put(6.2,7.2){$HN$}
\put(6,5){$\bullet$} \put(6.2,5.2){$\Stilde, \etilde, e', \Vtilde$}
\put(6,3){$\bullet$} \put(6.2,3.2){$S,\etilde, e, \Vtilde, V$} 
\put(6,1){$\bullet$}   \put(6.2,1.2){$N,\etilde,e, \Vtilde, V$}
%\put(6,-0.7){$\bullet$}   \put(6.2,-0.5){$1$}
\put(6.05,9.05){\line(0,-1){8}}

\put(3,6){$\bullet$} \put(1.7, 6.2){$H, \ctilde, c$}
\put(3,4){$\bullet$} \put(1.7, 4.2){$\Ttilde, \ftilde, f'$}
\put(3,2){$\bullet$} \put(1.3, 2.2){$T, \ftilde, f, U$}
\put(3,0){$\bullet$} \put(1.2, 0.2){$M, \ftilde, f, U$}
\put(3.05,6.05){\line(0,-1){6}}

\put(6.05,7.05){\line(-3,-1){3}}
\put(6.05,5.05){\line(-3,-1){3}}
\put(6.05,3.05){\line(-3,-1){3}}
\put(6.05,1.05){\line(-3,-1){3}}

\put(0,3){$\bullet$} \put(-0.4, 3.2){$R$}
\put(0,1){$\bullet$} \put(-0.4, 1.2){$Q$}
\put(0,-1){$\bullet$} \put(-0.7, -0.8){$\{1\}$}
\put(0.05,3.05){\line(0,-1){4}}

\put(3.05,4.05){\line(-3,-1){3}}
\put(3.05,2.05){\line(-3,-1){3}}
\put(3.05,0.05){\line(-3,-1){3}}

\end{picture}}
\end{center}

\bigskip

\begin{lemma}\label{lem isomofcomplexes}
For every $r\in R$, one has an isomorphism 
\begin{equation*}
   F\Stilde e\otimes_{FS} \left(Y_V\otimes_F FQ\otimes_F U^\circ \right)\otimes_{FT} fF\Ttilde
   \cong
   F\Stilde \lexp{r}e\otimes_{FS}\left(\lexp{r}{Y_V}\otimes_F FQ\otimes_F \lexp{r}{(U^\circ)} \right)\otimes_{FT} \lexp{r}fF\Ttilde
\end{equation*}
of chain complexes of $(F\Stilde e', F\Ttilde f')$-bimodules.
\end{lemma}

\begin{proof}
For any $M\in\lmod{FS}$, consider the map
\begin{equation*}
   F\Stilde e\otimes_{FS} \left(M\otimes_F FQ\otimes_F U^\circ \right)\otimes_{FT} fF\Ttilde
   \to
   F\Stilde \lexp{r}e\otimes_{FS}\left(\lexp{r}{M}\otimes_F FQ\otimes_F \lexp{r}{(U^\circ)} \right)\otimes_{FT} \lexp{r}fF\Ttilde\,,
\end{equation*}
mapping $a\otimes(y\otimes q\otimes u)\otimes b$ to $ar^{-1}\otimes (y\otimes rqr^{-1}\otimes u)\otimes rb$. It is straightforward to check that it is well-defined, an isomorphism of $(F\Stilde e', F\Ttilde f')$-bimodules, and functorial in $M$. Thus, it yields the desired isomorphism of chain complexes.
\end{proof}

\begin{nothing}
For the rest of the paper we assume that there exists $W\in \lmod{\FF_p Q}$ such that 
\begin{equation}\label{eqn hypo B}
   \text{{\em $\Res^S_Q V\cong F\otimes _{\FF_p} W$ and that $W$ has an endosplit $p$-permutation resolution $X_W$.}}
\end{equation}
Then the chain complex $F\otimes_{\FF_p} X_W$ is an endosplit $p$-permutation resolution of $\Res^S_Q V$ and we assume from now on that $X_V= F\otimes_{\FF_p }X_W$.

Note that that if $Q$ is abelian then (\ref{eqn hypo B}) is satisfied. In fact, every indecomposable endopermutation module for an abelian $p$-group is a direct summand of tensor products of inflations of Heller translates of the trivial module of quotient groups (see~\cite{Dade78} or \cite{Thevenaz2007}),  and every indecomposable endopermutation module (over any base field) is absolutely indecomposable (see Theorem~6.6 in the first paper \cite{Dade78}). It follows that $\Res^S_Q(V)$ has an $\FF_p$-form $W\in\lmod{\FF_pQ}$. Moreover, $W$ has an endosplit $p$-permutation resolution $X_W$ (see \cite[Theorem~7.2]{Rickard1996} whose proof is still valid over $\FF_p$).
\end{nothing}

\begin{proposition}\label{prop stableunderT}
Suppose that $R$ is abelian. For any $\stilde\in \Stilde$, one has an isomorphism $\lexp{\stilde}(\Ind^S_Q X_{W})\cong \Ind^S_Q X_{W}$ of complexes of $\FF_pS$-modules. In particular, for any $\stilde\in \Stilde$, one has $\lexp{\stilde}(\Ind^S_Q\Res^S_Q V)\cong \Ind^S_Q\Res^S_Q V$ as $FS$-modules. 
\end{proposition}

\begin{proof}
The complex $\Ind^S_Q X_{W}$ is isomorphic to a complex whose terms are direct sums of permutation $\FF_pS$-modules of the form $\FF_p[S/Q_0]$ where $Q_0\le Q$ and whose differentials are $\FF_p$-linear combination of maps of the form $f_t\colon\FF_p[S/Q_1]\to\FF_p[S/Q_2]$, $Q_1\mapsto \sum_{sQ_2\in Q_1tQ_2} sQ_2$, for some $t\in Q$, with $Q_1,Q_2\le Q$. Let $\stilde\in \Stilde=RN$ and write $\stilde=rn$ for some $r\in R$ and $n\in N$. For any $Q_0\le Q$, one has an isomorphism $\FF_p[S/Q_0]\to \lexp{\stilde\,}{\FF_p[S/Q_0]}$ of $\FF_pS$-modules given by $sQ_0\mapsto sn^{-1}Q_0$. Moreover, a quick computation shows that this isomorphism commutes with the above maps $f_t$, since $R$ is abelian. Therefore we have $\lexp{\stilde}(\Ind^S_Q X_{W})\cong \Ind^S_Q X_{W}$. For the last assertion note that this also implies that $\lexp{\stilde}(F\otimes_{\FF_p}\Ind^S_Q X_{W})\cong F\otimes_{\FF_p}\Ind^S_Q X_{W}$. Since the module $\Ind^S_Q\Res^S_Q V$ is the homology of the complex $F\otimes_{\FF_p}\Ind^S_Q X_{W}$ the result follows.
\end{proof}

\begin{lemma}\label{lem stabilizers}
Suppose that $R$ is abelian. Then one has $\stab_{\Gamma}(Z)=\stab_{\Gamma}(\omega)=p_1(X) =\stab_{\Gamma}(e')=\stab_{\Gamma}(f')=p_1(Y)$.
\end{lemma}

\begin{proof}
Note that $p_1(X) =\stab_{\Gamma}(e')=\stab_{\Gamma}(f')=p_1(Y)$ hold by (\ref{eqn p_1(X)=p_1(Y)}).

Since the complex $Z$ induces a splendid Rickard equivalence between $F\Stilde e'$ and $F\Ttilde f'$, the inclusion $\stab_{\Gamma}(Z)\le \stab_\Gamma(e')$ is immediate. Thus, $\stab_\Gamma(Z)\le p_1(X)$. 
Conversely, if $\sigma\in p_1(X)$, then $\lexp{\sigma}{e}=\lexp{\stilde}{e}$ for some $\stilde\in \Stilde$. Write $\stilde=rn$ for some $r\in R$ and $n\in N$ and note that we have $\lexp{\sigma}{e}=\lexp{r}{e}$.  This implies that $\lexp{\sigma}{V}\cong \lexp{r}{V}$ as $FS$-modules. By Proposition \ref{prop stableunderT}, we have $\lexp{\sigma}(F\otimes_{\FF_p}\Ind^S_Q X_{W})\cong F\otimes_{\FF_p} \Ind^S_Q X_{W}\cong \lexp{r}(F\otimes_{\FF_p} \Ind^S_Q X_{W})$ as complexes of $FS$-modules and $\lexp{\sigma}(\Ind^S_Q\Res^S_Q V)\cong \Ind^S_Q\Res^S_Q V\cong \lexp{r}(\Ind^S_Q\Res^S_Q V)$ as $FS$-modules. Therefore Lemma \ref{diagram} implies that $\lexp{\sigma}{Y_V} \cong \lexp{r}{Y_V}$ as complexes of $FS$-modules, as $Y_V$ was chosen to be contractible-free, see \ref{noth Rickard}.
Since the idempotents $e$ of $FS$ and $f$ of $FT$ are Brauer correspondents, also $\lexp{r}e$ and $\lexp{r}f$ are Brauer correspondents. Since the Galois action commutes with the Brauer correspondence, $\lexp{\sigma}{e}=\lexp{r}{e}$ implies $\lexp{\sigma}f=\lexp{r}f$. Therefore we have $\lexp{\sigma}U \cong \lexp{r}U$ as $FT$-modules.
By Lemma \ref{lem isomofcomplexes}, we obtain
\begin{align*}
   \lexp{\sigma}Z & \cong\lexp{\sigma}{\left(F\Stilde e\otimes_{FS}
         \left(Y_V\otimes_F FQ\otimes_F U^\circ \right)\otimes_{FT} fF\Ttilde\right)}\\
   & \cong F\Stilde \lexp{\sigma}e\otimes_{FS} \left(\lexp{\sigma}Y_V\otimes_F FQ\otimes_F \lexp{\sigma}U^\circ \right) 
            \otimes_{FT} \lexp{\sigma}fF\Ttilde\\
   & \cong F\Stilde \lexp{r}e\otimes_{FS} \left(\lexp{r}Y_V\otimes_F FQ\otimes_F \lexp{r}U^\circ \right)
               \otimes_{FT} \lexp{r}fF\Ttilde\\
   & \cong F\Stilde e\otimes_{FS} \left(Y_V\otimes_F FQ\otimes_F U^\circ \right)\otimes_{FT} fF\Ttilde \cong Z\,.
\end{align*}
This proves that $\stab_{\Gamma}(Z)=p_1(X)$. 

Since $\omega$ is a $p$-permutation equivalence between $F\Stilde e'$ and $F\Ttilde f'$, the inclusion $\stab_{\Gamma}(\omega)\le \stab_\Gamma(e') = p_1(X)$ is clear. The inclusion $\stab_{\Gamma}(Z)\le\stab_{\Gamma}(\omega)$ is immediate, and the proof is complete.
\end{proof}

\begin{corollary}\label{cor p-perm equiv}
Suppose that $R$ is abelian.

\smallskip
{\rm (a)} There exists a $p$-permutation equivalence between $\FF_p\Stilde \etilde$ and $\FF_p\Ttilde \ftilde$. 

\smallskip
{\rm (b)} There exists a $p$-permutation equivalence between $\FF_pG\btilde$ and $\FF_pH\ctilde$.
\end{corollary}

\begin{proof}
(a) By Lemma \ref{lem stabilizers} we have $\stab_{\Gamma}(\omega)=\stab_{\Gamma}(e')=\stab_{\Gamma}(f')$. Hence by Lemma \ref{lem descending} there exists a $p$-permutation equivalence between $\FF_p\Stilde\etilde$ and $\FF_p\Ttilde\ftilde$. 

\smallskip
(b) The $p$-permutation bimodule $\FF_pG\etilde$ induces a Morita equivalence, hence a $p$-permutation equivalence, between $\FF_pG\btilde$ and $\FF_p\Stilde\etilde$. Similarly, the bimodule $\FF_pH\ftilde$ induces a $p$-permutation equivalence between $\FF_pH\ctilde$ and $\FF_p\Ttilde\ftilde$. The result follows now from Part~(a). 
\end{proof}

\begin{Remark}\label{rem no descent theory for Z}
If one had a descent result for splendid Rickard equivalences analogous to Lemma~\ref{lem descending}, one would also obtain a splendid Rickard equivalence between $\FF_p G\btilde$ and $\FF_p H\ctilde$ in the above corollary, because Lemma~\ref{lem stabilizers} includes $\stab_\Gamma(Z)$, while in the proof of the above corollary we only used the statement about $\stab_\Gamma(\omega)$. In order to prove such a descent result one would need a descent result for homomorphisms between $p$-permutation modules.

Moreover, the approach in the proof of Proposition~\ref{prop Rickard equivalence over Fp} does not work, since the first Fong reduction only gives an equivalence between $\FF_pG\btilde$ and $\FF_p\Stilde\etilde=\FF_p\Stilde{\tilde{e'}}$, and not between $\FF_p G\btilde$ and $\FF_p S\etilde$. In order to apply the descent result from \cite{KL2018}, one would first need to descend the chain complex $Z$ from \ref{noth Z}(d) to $\FF_p[e']$. But we could not modify the approach from the proof of Proposition~\ref{prop Rickard equivalence over Fp} to descend to $\FF_p[e']$, unless $e=e'$ which is equivalent to $Q=R$ and to $\FF_p[b]=\FF_p[e]$.
\end{Remark}

%%%%%%%%%%%% REFERENCES %%%%%%%%%%%%%%%%%%%%%%%%%%%%%%%%%%%%%%%%%

\end{document}